\documentclass[12pt]{amsart}
\usepackage{amsmath,amsfonts,amssymb,amscd}
\usepackage{latexsym,epsfig}
\usepackage[latin1]{inputenc}
\usepackage[T1]{fontenc}
\usepackage{enumerate}
\usepackage{amsbsy}
\usepackage[matrix,arrow]{xy}

\newcommand{\N}{\ensuremath{\mathbb N}}
\newcommand{\Z}{\ensuremath{\mathbb Z}}
\newcommand{\F}{\ensuremath{\mathbb F}}
\newcommand{\T}{\ensuremath{\mathbb T}}
\newcommand{\Sf}{\ensuremath{\mathbb S}}

\newtheorem{thm}{Theorem}
\newtheorem{lem}[thm]{Lemma}
\newtheorem{prop}[thm]{Proposition}
\newtheorem{cor}[thm]{Corollary}

\newcommand{\reth}[1]{Theorem~\protect\ref{th:#1}}
\newcommand{\relem}[1]{Lemma~\protect\ref{lem:#1}}

\begin{document}

\title{On residual properties of  pure braid groups of closed surfaces}

\author[Bardakov]{Valerij G. Bardakov}
\address{Sobolev Institute of Mathematics, 630090 Novosibirsk, Russia}
\email{bardakov@math.nsc.ru}

\author[Bellingeri]{Paolo Bellingeri}
\address{Univ. Milano Bicocca, Dipartimento Matematica e Applicazioni, 20126 Milan, Italy}
\email{paolo.bellingeri@unimib.it}



\begin{abstract}
We prove that pure braid groups of closed surface are almost-direct products
of residually torsion free nilpotent groups and hence residually torsion free nilpotent.
As a Corollary, we prove also that braid groups on $2$ strands of closed surfaces are residually nilpotent.
\end{abstract}

\maketitle





\section{Introduction}

A group $G$ is said to be
\emph{residually torsion free nilpotent} if for any (non-trivial) element $x
\in G$, there exists a homomorphism $\phi: G \to H$ such that   $H$
is  torsion free  nilpotent and  $\phi(x) \not=1$.

Let $A, C$ be two groups. If   $C$ acts on $A$
by
automorphisms, the semi-direct product $A \rtimes C$ is said to be
\emph{almost-direct} if the action of $C$ on the abelianization
of $A$ is trivial. An example
 of  almost-direct product of free groups is given by Artin pure braid group $P_n$.
Such decomposition implies that $P_n$ is residually torsion-free nilpotent (see Section 2).

The structure of almost-direct product turns out to be also a powerful tool also in the determination of algebras
related to lower central series (see for instance \cite{CCP}) and more generally in the study of finite type invariants.
The decomposition  of $P_n$ as almost-direct product
of free groups was used in~\cite{pap}
in order to construct an universal finite type invariant
for braids with integers coefficients.

Let  $\Sigma$ be an oriented  surface and let
$\F_n(\Sigma)=\{(x_1, \dots x_n) \in \Sigma^n \, | \, x_i \not= x_j$
for $i\not=j \}$. The fundamental group of $\F_n(\Sigma)$ is
called \emph{pure braid group} on $n$ strands
of $\Sigma$ and it is usually denoted by $P_n(\Sigma)$. When
$\Sigma$ is the disk $D^2$ we obtain a group which  is isomorphic to Artin
pure braid group $P_n$.

The symmetric group $S_n$ acts on $\F_n(\Sigma)$ by permutation of coordinates
and the fundamental group of the orbit space  $\F_n(\Sigma)/S_n$ is called \emph{braid
 group}
 on  $n$ strands of $\Sigma$ and it is denoted by $B_n(\Sigma)$.
For $n=1$ we have that  $P_1(\Sigma)=B_1(\Sigma)=\pi_1(\Sigma)$ and $B_n(D^2)$ is isomorphic to Artin braid group
$B_n$.

When $\Sigma$ is an oriented  surface of positive genus,  $P_n$ embeds naturally
into $P_n(\Sigma)$. In~\cite{GP} Gonz\'alez-Meneses and Paris proved that
the  normal closure of the classical pure braid group $P_n$ in $P_n(\Sigma)$
is an almost-direct product of (infinitely generated) free groups. Adapting the approach of Papadima,
they constructed a universal finite type invariant for surface braids.

At our knowledge, it is not known if pure braid groups of surfaces (different from the disk)
can be decomposed as almost-direct products of (residually) free groups
(see for instance~\cite{B,GG} for more details on this subject).
In  \cite{BGG} the first author proved that  pure braid
groups of the torus and of surfaces with boundary components are
residually torsion-free nilpotent by showing that
they may be realised as  subgroups of the Torelli group of a surface
of higher genus (see also the end of Section~\ref{mainsection}).

In this paper we complete the study of    lower central series   and  related
residual properties of (pure) braid groups of surfaces begun in \cite{BGG}, proving that
pure braid groups of closed surface are almost-direct products
of residually torsion free nilpotent groups and hence residually torsion free nilpotent
(Theorem~\ref{th:main}).
As a Corollary, we prove also that braid groups on $2$ strands of closed surfaces are residually nilpotent
(Corollary~\ref{b2}).

The fact that a group is residually torsion-free nilpotent has several
 consequences, notably that the group is bi-orderable~\cite{MR}
and residually $p$-finite~\cite{Gr}. Therefore it follows from Theorem~\ref{th:main}
that pure braid groups of closed oriented surfaces
 are  bi-orderable
and residually $p$-finite; the first result was earlier proved in~\cite{Go} and the second is also a consequence of  Theorem~1.2 in~\cite{Paris} (see  Section~\ref{mainsection}).

\vspace*{5pt}

\noindent{\bf Acknowledgments.} The research of the first author
has been supported by the University of Nantes. The first author would like to thank the
members of the Department of Mathematics of the  University of Nantes
for their kind hospitality.

\section{Residual properties, almost-direct products and group presentations for pure braid
groups on closed surfaces}

Let us begin with few definitions.

The lower central series of a group $G$ is the filtration
 $\Gamma_1(G)=G \supseteq \Gamma_2(G) \supseteq \ldots$, where
$\Gamma_i(G)=[\Gamma_{i-1}(G), G]$. The \emph{rational lower central series}
of $G$ is the   filtration $D_1(G) \supseteq  D_2(G)
\supseteq \ldots$  obtained setting $D_1(G)=G$, and for $i\geq 2$,
defining $D_i(G)=\{ \, x \in G \, | \, x^n \in \Gamma_{i}(G) $ for some
$n\in \N \setminus \{ 0 \} \, \}$.

Let $\mathcal{F}\mathcal{P}$ be the family of groups
having the group-theoretic property $\mathcal{P}$. A group $G$ is said to be
\emph{residually $\mathcal{P}$} if for any element $x
\in G \setminus 1$, there exists a  homomorphism of $G$ into some group in $\mathcal{F}\mathcal{P}$
 taking  $x$ in a nontrivial element.

A group $G$ is residually nilpotent if and only if
$\bigcap_{i \ge 1}\Gamma_i(G)=\{ 1\}$. On the other hand, a group $G$
is residually torsion-free nilpotent if and only if
$\bigcap_{i \ge 1}D_i(G)=\{ 1\}$.

\begin{prop}{\bf (\cite{bb1,FR2})}\label{corrn}
Let $A, C$ be two groups such that $C$ acts on $A$ by automorphisms. If $A \rtimes C$ is an \emph{almost-direct product} then
$\Gamma_m(A \rtimes C)=\Gamma_m(A) \rtimes \Gamma_m(C)$ and
$D_m(A \rtimes C)=D_m(A) \rtimes D_m(C).$
\end{prop}

\begin{cor}\label{corrn2}
The almost-direct product of two residually nilpotent (torsion free) groups is
residually nilpotent (torsion free).
\end{cor}

The pure braid group $P_n$ is an almost-direct product
of free groups (\cite{FR}).
Since free groups are residually torsion-free nilpotent~\cite{F},
it follows from Corollary~\ref{corrn2} that pure braid groups are
residually torsion-free nilpotent~(see also \cite{FR2}).





\section{Presentations for pure braid groups on surfaces}
Let $\Sigma_g$ be an oriented closed surface of genus $g$.
Let $\mathcal{X}=\{x_1, \dots,
x_n\}$ be a set of $n$ distinct points (\emph{punctures}) in the
interior of $\Sigma_g$. A \emph{pure geometric braid} on $\Sigma_g$ based at $\mathcal{X}$
is a collection $(\psi_1, \dots, \psi_n)$ of $n$ disjoint paths
(called \emph{strands}) on $\Sigma_g \times [0, 1]$
 which run
monotonically with $t \in [0, 1]$ and such that
$\psi_i(0)=(x_i, 0)$ and $\psi_i(1) =(x_i, 1)$. Two pure braids are considered to be equivalent if they
are isotopic relatively to the base points. The usual product of paths defines a group structure on
the equivalence classes of braids. This group, which is isomorphic to
$P_n(\Sigma_g)$, does not depend on the choice of $\mathcal{X}$.

We recall a group presentation for pure braid groups of oriented closed surfaces~\cite{B}.
In the following we set  $[a,b] =a^{-1} b^{-1} a b$, $a^b=b^{-1} a b$
and ${}^{b}a=b a b^{-1}$ and  we use  the convention that 
$W = \prod_{i=m}^{n} f(i)$ with $n < m$ implies that $W =1$. 

\begin{thm}\textbf{(\cite{B})} \label{th:pp}
Let $g \ge 1$ and $n \ge 2$. The group $P_n(\Sigma_g)$ admits the following
presentation:

Generators: $\{A_{i,j}\; | \;1 \le i \le 2g+ n -1,  2g +1\le j \le 2g + n, i<j \}.$

Relations:
\begin{eqnarray*}
 &\text{(PR1)}&  A_{i,j}^{-1}  A_{r,s} A_{i,j} = A_{r,s} \;  \; \mbox{if} \, \,(i<j<r<s)  \;   \mbox{or} \,
(r+1<i<j<s),\\
 &     &  \mbox{or} \, (i=r+1<j<s \, \, \mbox{for even} \, \, r<2g  \, \, \mbox{or}   \, \, r>2g  ) \,; \\
 &\text{(PR2)}&  A_{i,j}^{-1}  A_{j,s} A_{i,j} = A_{i,s}  A_{j,s} A_{i,s}^{-1} \;  \; \mbox{if} \, \,
(i<j<s)\,;\\
 &\text{(PR3)}&  A_{i,j}^{-1}  A_{i,s} A_{i,j} = A_{i,s} A_{j,s} A_{i,s}  A_{j,s}^{-1} A_{i,s}^{-1} \; \;
\mbox{if} \, \, (i<j<s)\,; \\
 &\text{(PR4)}&
A_{i,j}^{-1}A_{r,s}A_{i,j}=A_{i,s}A_{j,s}A_{i,s}^{-1}A_{j,s}^{-1}A_{r,s}A_{j,s}A_{i,s}A_{j,s}^{-1}A_{i,s}^{-1}
\\
 &     &   \mbox{if} \, \,(i+1<r<j<s)  \;   \mbox{or} \\
 &     &   \, \, (i+1=r<j<s  \; \mbox{for odd }  \, \, r<2g  \, \, \mbox{or}   \, \, r>2g ) \,;
 \end{eqnarray*}
 \begin{eqnarray*}
 &\text{(ER1)}&  A_{r+1,j}^{-1} A_{r,s} A_{r+1,j}=A_{r,s} A_{r+1,s} A_{j,s}^{-1} A_{r+1,s}^{-1}  \\
 &     &     \mbox{if} \, \, j<s ,    \,r \, \mbox{odd and}\, \,r<2g                 \,               ; \\
 &\text{(ER2)}&   A_{r-1,j}^{-1}  A_{r,s} A_{r-1,j}=  A_{r-1,s} A_{j,s} A_{r-1,s}^{-1}  A_{r,s} A_{j,s}
A_{r-1,s} A_{j,s}^{-1}A_{r-1,s}^{-1} \\
 &    &   \mbox{if} \, \, j<s ,     \,r \, \mbox{even and}\, \,r<2g   \,               . \\
 \end{eqnarray*}
\begin{eqnarray*}
 &\text{(TR)}&  (\prod_{i=1}^g[A_{2i-1, 2g+k}^{-1},A_{2i, 2g+k}])^{-1} =
\prod_{l=2g+1}^{2g+k-1} A_{l,2g+k} \prod_{j=2g+k+1}^{2g+n} A_{2g+k,j}  \\
& & \quad  (k=1, \dots, n)    \,.
\end{eqnarray*}
\end{thm}

\begin{figure}[h]
\begin{center}
\psfig{figure=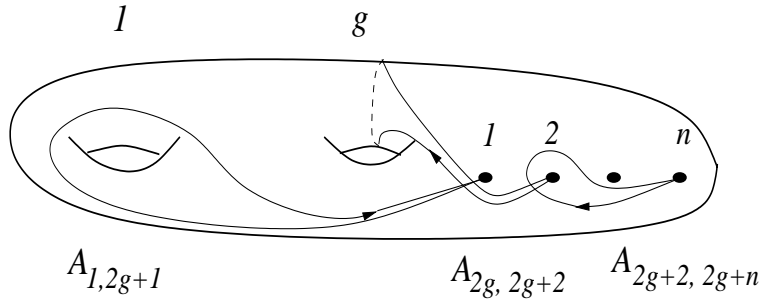,width=10cm,height=4cm}
\caption{\label{generateur}Projection of representatives of the generators $A_{i,j}$. We
represent $A_{i,j}$ by its only non-trivial strand.}
\end{center}
\end{figure}

As a representative of the generator $A_{i,j}$, we may take a
geometric braid whose only non-trivial (non-vertical) strand is the $(j-2g)$th
one.  In Figure~\ref{generateur}, we illustrate the projection of such braids on the surface $\Sigma_{g}$
(see also Figure~8 of~\cite{B}). Some misprints in Relations~(ER1) and~(ER2) of Theorem~5.1 of~\cite{B} have
been corrected. Remark also that in~\cite{B} was used the convention 
 $[a,b] =a b a^{-1} b^{-1}$. 

As remarked in the proof of  Theorem~5.1 of~\cite{B}, using above list of relations we can write
any element of type $A_{i,j} A_{p,q} A_{i,j}^{-1}$ for  $2g+1 \le j<q \le 2g+n$ as a word on the generators
$ A_{1,q}, \dots,   A_{q-1,q}$. Therefore starting from the group presentation in~\reth{pp}, we can obtain the following group presentation
for $P_n(\Sigma_g)$, where, in respect to~\reth{pp},
we set $A_{2r-1,2g+j}$=$a_{j,r}$ and $A_{2r,2g+j}=b_{j,r}$ for $r=1, \ldots, g$ and $j=1, \ldots, n$
and $A_{2g+i,2g+j}= \tau_{i,j}$ for $1 \le i < j \le  n$. 

\begin{prop}\label{prop:pp}
Let $g \ge 1$ and $n \ge 2$. The group $P_n(\Sigma_g)$ admits the following
presentation:

Generators: $\{a_{j,k}, \, b_{j,k},\, \tau_{p,q} | \;1 \le k \le g,  1 \le k  \le n, \,   1 \le p<q  \le n \}.$

Relations:
\begin{eqnarray*}
&\text{(I-a)}&  {}^{a_{j,k}}  c_{l,m}  = c_{l,m}  \qquad  \mbox{for} \qquad c_{l,m}=a_{l,m}, b_{l,m} \;
 \mbox{if} \, \, (m<k)\,;\\
&\text{(I-b)}&  {}^{b_{j,k}}  c_{l,m}  = c_{l,m}  \qquad  \mbox{for} \qquad c_{l,m}=a_{l,m}, b_{l,m} \;
 \mbox{if} \, \, (m<k)\,;\\
&\text{(I-$\tau$1)}& {}^{\tau_{s,j}} c_{l,m} = c_{l,m}  \qquad  \mbox{for} \qquad c_{l,m}=a_{l,m}, b_{l,m} \,;\\
&\text{(I-$\tau$2)}& {}^{\tau_{p,j}} \tau_{s,\,l } = \tau_{s,\,l } \, ;\\
&\text{(II-a)}&  {}^{a_{j,k}} a_{l,k}  = a_{l,k}^{\tau_{j,\,l }}  ;\\
&\text{(II-b)}&  {}^{b_{j,k}} b_{l,k}  = b_{l,k}^{\tau_{j,\,l }}  ;\\
&\text{(II-$\tau$)}&  {}^{\tau_{s,j}} \tau_{s,\,l }  = \tau_{s,\,l }^{\tau_{j,\,l }}  ;\\
&\text{(III-a1)}&  {}^{a_{j,k}} c_{l,m}  = {}^{[\tau_{j,\,l }, a_{l,k}]} c_{l,m}  \qquad \mbox{for} \qquad c_{l,m}=a_{l,m}, b_{l,m}
 \mbox{if} \, \, (k<m)\,;\\
&\text{(III-a2)}&    {}^{a_{j,k}} \tau_{s,\,l }= {}^{[\tau_{j,\,l }, a_{l,k}]} \tau_{s,\,l };\\
&\text{(III-b1)}&  {}^{b_{j,k}} c_{l,m}  = {}^{[\tau_{j,\,l }, b_{l,k}]} c_{l,m}  \qquad \mbox{for} \qquad c_{l,m}=a_{l,m}, b_{l,m}
 \mbox{if} \, \, (k<m)\,;\\
&\text{(III-b2)}&    {}^{b_{j,k}} \tau_{s,\,l }= {}^{[\tau_{j,\,l }, b_{l,k}]} \tau_{s,\,l };\\
&\text{(III-$\tau$)}&  {}^{\tau_{s,j}} \tau_{p,\,l }  = {}^{[\tau_{j,\,l }, \tau_{s,\,l }]} \tau_{p,\,l } \,;\\
&\text{(IV-a)}&  {}^{a_{j,k}} \tau_{r,\,l }  = \tau_{r,\,l };\\
&\text{(IV-b)}&  {}^{b_{j,k}} \tau_{r,\,l }  = \tau_{r,\,l };\\
&\text{(IV-$\tau$)}&  {}^{\tau_{p,j}} \tau_{r,\,l }  = \tau_{r,\,l };\\
&\text{(V-a)}&  {}^{a_{j,k}} \tau_{j,\,l }  =  [\tau_{j,\,l }, a_{l,k}]  \tau_{j,\,l }  ;\\
&\text{(V-b)}&  {}^{b_{j,k}} \tau_{j,\,l }  =  [\tau_{j,\,l }, b_{l,k}]  \tau_{j,\,l }  ;\\
&\text{(V-$\tau$)}&  {}^{\tau_{s,j}} \tau_{j,\,l }  =  [\tau_{j,\,l }, \tau_{s,\,l }]  \tau_{j,\,l }  ;\\
&\text{(ER1)}&  {}^{a_{j,k}} b_{l,k}  = \tau_{j,\,l }^{-1} b_{l,k} [a_{l,k}, \tau_{j,\,l }]  ;\\
&\text{(ER2)}&  {}^{b_{j,k}} a_{l,k}  = a_{l,k} \tau_{j,\,l }  ;\\
&\text{(TR)}& \prod_{i=1}^g[a_{l,i}^{-1}, b_{l,i}]=
(\prod_{w=1}^{l-1} \tau_{w,\,l } \prod_{d=l+1}^{n} \tau_{l,d})^{-1};\\
\end{eqnarray*}
where $1 \le m,k \le g$ and  $1\le s< p < j < r < l < q  \le n$ (with $s, p, q, r$ and $l$ possibly absent).
\end{prop}

Remark also that Proposition~\ref{prop:pp} can be proven directly. In fact drawing corresponding braids,
one can verify that above relations hold in $P_n(\Sigma_g)$. To prove that they form a complete set of relations
it is sufficient to repeat the arguments in the proof of  Theorem~5.1 of~\cite{B}.





\section{The structure of pure braid groups of closed surfaces}\label{mainsection}

Let $p:P_n(\Sigma_g) \to \pi_1(\Sigma_g)$ be the map which forgets all strands except the first one.
This map is induced from the forgetting map at the level of corresponding configuration spaces and
 $\ker p$ is isomorphic to $P_{n-1}(\Sigma_{g,1})$, the pure braid group on $n-1$ strands of the oriented surface
$\Sigma_{g,1}$ of genus $g$ with one boundary component (see for instance~\cite{GG}).

In the following we provide an  algebraic section $s$ for $p:P_n(\Sigma_g) \to \pi_1(\Sigma_g)$, where
$ P_n(\Sigma_g)$ has the  presentation given in Proposition~\ref{prop:pp} and
 we show  that $s$ induces a structure of almost-direct product on
$P_n(\Sigma_g)$.

In~Theorem~1 of~\cite{GG} it was shown that $p$ admits a geometric  section (i.e. induced from a section
on the geometric level of corresponding  configuration spaces).  In~\cite{GG} it was also given an algebraic definition of such section using another group
presentation for $P_n(\Sigma_g)$ (provided in Corollary~8 of~\cite{GG}).

Before stating the main result of the paper we need a preliminary Lemma.
In the following, we will set  $T_{l,\,q} = \prod_{d=l+1}^{q} \tau_{l,\,d}$ for $1\le   l<q  \le n$.

\begin{lem}\label{lem:1}
The following identities:
\begin{eqnarray}
&(\prod_{d=l}^{q} c_{d,k}) T_{l,\,q}=T_{l,\,q} (\prod_{d=l}^{q} c_{d,k})& \qquad \text{for} \qquad
c_{d,k} = a_{d,k} , \, b_{d,k} \, ;\\
&{}^{b_{l,k}} (\prod_{d=l+1}^{q} a_{d,k})= (\prod_{d=l+1}^{q} a_{d,k}) T_{l,\,q}& \, ,
\end{eqnarray}
hold in $P_n(\Sigma_g)$  for $1\le   l<q  \le n$ and $1 \le k \le g$.
\end{lem}
\begin{proof}
We provide an algebraical verification of these identities  that  can be also verified drawing corresponding braids.

We prove the first identity, in the case of $c_{d,k} = a_{d,k}$ (the case of $c_{d,k} = b_{d,k}$ is analogous).
First we remark that for $1\le   i<l<j  \le n$ and $1 \le k \le g$,

\begin{eqnarray*}
 a_{l,k} \underline{\tau_{i,\,\,l } \tau_{i,\,j} \tau_{i,\,\,l }^{-1}} a_{l,k}^{-1}
= 
  \underline{a_{l,k} \tau_{l,\,j}^{-1} \tau_{i,\,j} \tau_{l,\,j} a_{l,k}^{-1}} =
 \tau_{l,\,j}^{-1} \tau_{i,\,j} \tau_{l,\,j}  \, ,
\end{eqnarray*}

respectively because of relation II-$\tau$, relation III-a2 and relation V-a.
Therefore applying once more  relation II-$\tau$ we obtain that 

\begin{eqnarray*}
(A) \quad (\tau_{i,\,\,l }^{-1} a_{l,k} \tau_{i,\,\,l }) \tau_{i,\,j} (\tau_{i,\,\,l }^{-1} a_{l,k}^{-1}
\tau_{i,\,\,l })= \underline{\tau_{i,\,\,l }^{-1} \tau_{l,\,j}^{-1} \tau_{i,\,j} \tau_{l,\,j} \tau_{i,\,\,l }} =\tau_{i,\,j}
\end{eqnarray*}

holds in $P_n(\Sigma_g)$ for $1\le   i<l<j  \le n$ and $1 \le k \le g$.
Now we claim that for $1\le   l<q  \le n$ and $1 \le k \le g$,

\begin{eqnarray*}
(B) \quad (\prod_{d=l}^{q} a_{d,k}) \tau_{l,\,q}=\tau_{l,\,q} (\prod_{d=l}^{q} a_{d,k}) \, .
\end{eqnarray*}

In fact,

$$ (\prod_{d=l}^{q} a_{d,k}) \tau_{l,\,q} (\prod_{d=l}^{q} a_{d,k})^{-1}= (\prod_{d=l+1}^{q-1}{}^{a_{l,k}}a_{d,k}) \, \cdot \,
{}^{a_{l,k}}(a_{q,k} \tau_{l,\,q} a_{q,k}^{-1})  \, \cdot \, (\prod_{d=l+1}^{q-1} {}^{a_{l,k}}a_{d,k})^{-1}=
$$

and by relation II-a,

$$
=(\prod_{d=l+1}^{q-1} a_{d,k}^{\tau_{l,d}})\, \cdot \,
\underline{{}^{a_{l,k}}(a_{q,k} \tau_{l,\,q} a_{q,k}^{-1})} \, \cdot \, (\prod_{d=l+1}^{q-1} a_{d,k}^{\tau_{l,d}})^{-1}=
$$

and therefore applying relations II-a and V-a,

$$
=(\prod_{d=l+1}^{q-1} a_{d,k}^{\tau_{l,d}})
 (\tau_{l,\,q})  (\prod_{d=l+1}^{q-1} a_{d,k}^{\tau_{l,d}})^{-1} \,.
$$

Applying   identity  (A) we obtain that   $(\tau_{l,d}^{-1} a_{d,k} \tau_{l,d})
 \tau_{l,\,q}  (\tau_{l,d}^{-1}  a_{d,k}^{-1} \tau_{l,d}) = \tau_{l,\,q}$ for $l+1 \le d <q \le n$
and therefore

$$
(\prod_{d=l+1}^{q-1} a_{d,k}^{\tau_{l,d}})
 (\tau_{l,\,q}) (\prod_{d=l+1}^{q-1} a_{d,k}^{\tau_{l,d}})^{-1} =\tau_{l,\,q}
$$

and the identity (B) is proved.
Thus,
because of relation I-$\tau$1  the following identity holds:

$$
(\prod_{d=l}^{q} a_{d,k}) T_{l,\,q}= a_{l,k} a_{l+1,k} \tau_{l,l+1} a_{l+2,k} \tau_{l,l+2} \cdots   a_{q,k} \tau_{l,q}
$$

Hence applying relation (B) recursively   we obtain that

$$
( \prod_{d=l}^{q} a_{d,k} ) T_{l,\,q}= \prod_{d=l+1}^{q} \tau_{l, \, d} \prod_{d=l}^{q} a_{d,k}
=T_{l,\,q} (\prod_{d=l}^{q} a_{d,k}) \, .
$$

The second identity is easier to verify. From relations ER2 and I-$\tau$1
 in Proposition~\ref{prop:pp} one obtains the following identities:
\begin{eqnarray*}
&{}^{b_{l,k}}&(\prod_{d=l+1}^{q} a_{d,k})=\prod_{d=l+1}^{q} ( a_{d,k} \tau_{l,d})=
a_{l+1,q} (\prod_{d=l+1}^{q-1} (\tau_{l,d} a_{d+1,k}))\tau_{l,q}= \\
&=& \cdots =\prod_{d=l+1}^{q} (a_{d,k}) T_{l,q} \, .
\end{eqnarray*}
\end{proof}

Let $\langle c_1, d_1, \ldots, c_g, d_g \, | \prod_{i=1}^g[c_{i}^{-1}, d_{i}]=1 \rangle$ be a group presentation for
$\pi_1(\Sigma_g)$; the morphism $p:P_n(\Sigma_g) \to \pi_1(\Sigma_g)$ can be defined algebraically as follows:
$p(a_{1,k})=c_k$ and  $p(b_{1,k})=d_k$ for $1 \le k \le g$ and $p(a_{j,k})=p(b_{j,k})=\tau_{p,q}=1$ elsewhere.

\begin{thm}\label{th:main}
The exact sequence
\begin{eqnarray}\label{eqn:1}
1 \to P_{n-1}(\Sigma_{g,1}) \to P_n(\Sigma_g) \to \pi_1(\Sigma_g) \to 1
\end{eqnarray}
splits and $P_n(\Sigma_g) \simeq P_{n-1}(\Sigma_{g,1}) \rtimes \pi_1(\Sigma_g)$ is almost-direct product of $P_{n-1}(\Sigma_{g,1})$  and $\pi_1(\Sigma_g)$.
\end{thm}

\begin{proof}
Let us define a set-section $s:\pi_1(\Sigma_g) \to  P_n(\Sigma_g)$ as follows:
 $s(c_k)=T_{1,n} a_{1,k} T_{1,n}^{-1}$ and $s(d_k)=T_{1,n} b_{1,k} T_{1,n}^{-1}$
for $1 \le k <g$, $s(c_g)= \prod_{d=1}^{n} a_{d,g} T_{1,\,n}$ and $s(d_g)= b_{1,g}$.
In order to prove that $s$ is a well-defined morphism  
it suffices to prove that  
$$
\prod_{i=1}^g[s(c_{i})^{-1}, s(d_{i})]=s(\prod_{i=1}^g[c_{i}^{-1},d_{i}])=1 \, .
$$

From \relem{1}  and relation TR in Proposition~\ref{prop:pp}  one deduces the following identities:

$$
\prod_{i=1}^g[s(c_{i}^{-1}), s(d_{i})]=T_{1,n} \prod_{i=1}^{g-1}[a_{1,i}^{-1}, b_{1,i}] T_{1,n}^{-1}
\prod_{d=1}^{n} a_{d,g} T_{1,\,n} b_{1,\,g}^{-1} (\prod_{d=1}^{n} a_{d,g} T_{1,\,n})^{-1} b_{1,\,g}=
$$

$$
=T_{1,n} \prod_{i=1}^{g-1}[a_{1,i}^{-1}, b_{1,i}]
\prod_{d=1}^{n} a_{d,g} b_{1,g}^{-1} (\prod_{d=1}^{n} a_{d,g} T_{1,\,n})^{-1} b_{1,g}=T_{1,n} \prod_{i=1}^{g-1}[a_{1,i}, b_{1,i}^{-1}] \times
$$
$$
\times a_{1,g} b_{1,g}^{-1} \cdot
{}^{b_{1,g}}(\prod_{d=2}^{n} a_{d,g})  \cdot (\prod_{d=2}^{n} a_{d,g} T_{1,\,n})^{-1} a_{1,g}^{-1} b_{1,g}=T_{1,n} \prod_{i=1}^{g}[a_{1,i}^{-1}, b_{1,i}] = 1 \,
$$

and therefore $s:\pi_1(\Sigma_g) \to  P_n(\Sigma_g)$ is a well-defined morphism and~(\ref{eqn:1}) splits.

Now, remark that from relation ER2 we deduce the following identities: 
$$
\tau_{w,\,l }=[a_{l,g},\,b_{w,g}^{-1}] \quad \mbox{for} \; 1 \le w <l \le n \; ;
$$
$$
\tau_{l,d}=[a_{d,g},\,b_{l,g}^{-1}] \quad \mbox{for} \; 1 \le l <d \le n \;.
$$

Therefore, from relations TR and ER2 in  Proposition~\ref{prop:pp} we obtain that the following relation
$$
\tau_{1,\,l }^{-1} = \prod_{w=2}^{l-1} \tau_{w,\,l } \prod_{d=l+1}^{n} \tau_{l,d}  \prod_{i=1}^g[a_{l,i}^{-1}, b_{l,i}]=
$$
$$
=\prod_{w=2}^{l-1} [a_{l,g},\, b_{w,g}^{-1}] \prod_{d=l+1}^{n} [a_{d,g},\, b_{l,g}^{-1}] \prod_{i=1}^g[a_{l,i}^{-1}, b_{l,i}] \,
$$
holds in $P_n(\Sigma_g)$ for $2\le  l \le n$ and then   $\tau_{1,\,l } \in \Gamma_2(\ker p)$ for $l=2, \ldots,n$.

Fixing $j=1$, relations in Proposition~\ref{prop:pp} provide  the action
by conjugacy of $a_{1,k}, b_{1,k}$ for $1\le k \le g$ on  the set $\mathcal{A}=
\{ a_{j,k}, b_{j,k}, \tau_{p,q}  \, | \, 1\le k \le g$, $2 \le j \le n$ and  $1 \le p<q \le n \}$
and using the fact that $\tau_{1,\,l } \in \Gamma_2(\ker p)$ for $l=2, \ldots,n$
one can easily check that

$${}^{a_{1,k}} h \equiv h \qquad mod  \; \Gamma_2(\ker p)$$

and

$${}^{b_{1,k}} h \equiv h \qquad mod  \; \Gamma_2(\ker p)$$

for $1\le k \le g$ and $h \in \mathcal{A}$. 

Hence the action of $\pi_1(\Sigma_g)$  on the abelianisation of $\ker p$ is trivial.
In fact, let $h \in \mathcal{A}$. 
It follows from previous congruences that

$$ {}^{s(c_k)}h ={}^{T_{1,n} a_{1,k} T_{1,n}^{-1}}h \equiv h \qquad mod \; \Gamma_2(\ker p)$$

and

$$ {}^{s(d_k)} h ={}^{T_{1,n} b_{1,k} T_{1,n}^{-1}} h \equiv h \qquad mod  \; \Gamma_2(\ker p) \, ,$$

for $1\le k \le g-1$ and $h \in \mathcal{A}$.
On the other hand since also $\prod_{d=2}^{n} a_{d,g} \in \ker p$, one derives that

$$ {}^{s(c_g)} h ={}^{\prod_{d=1}^{n} a_{d,g} T_{1,\,n}} h \equiv {}^{a_{1,g}} h \equiv h \qquad mod  \; \Gamma_2(\ker p) \,.$$

Finally $s(d_g)=b_{1,g}$ and therefore

$$ {}^{s(d_g)} h \equiv h \qquad mod \Gamma_2(\ker p) \, .$$

Since  $\mathcal{A}$
is  a complete set of generators for $\ker p$, the action of $\pi_1(\Sigma_g)$ is  trivial on the abelianisation of $\ker p$.
\end{proof}

\begin{cor}
The group $P_n(\Sigma_g)$ is residually torsion free nilpotent for $n\ge 1$ and $g>0$.
\end{cor}
\begin{proof}
The group $P_1(\Sigma_g)$ is isomorphic to $\pi_1(\Sigma_g)$ which is residually free and therefore
residually torsion free nilpotent. In the case $n>1$ the claim is a straightforward consequence of~Corollary~\ref{corrn2},
\reth{main} and the fact that $P_{n-1}(\Sigma_{g,1})$ is residually torsion free nilpotent for $n\ge 1$ and $g,p>0$~\cite{BGG}.
\end{proof}

We remark that was already proven in~\cite{BGG}  that $P_n(\T^2)$ is residually torsion free nilpotent and that
$P_n(\Sf^2)$ is residually nilpotent but not residually torsion free nilpotent.  The group $P_n(\Sigma_g)$
was proven to be bi-orderable in~\cite{Go}.

We recall that residual torsion free nilpotence implies the residual $p$-finiteness.
We recall that the Torelli group $\mathcal{T}(\Sigma_g)$ of the surface $\Sigma_g$
is defined as the kernel of the natural action of the mapping class group
of $\Sigma_g$ on $H_1(\Sigma_g)$. Let $\mathcal{P}$ a set of $n$ distinct point on $\Sigma_g$.
According to~\cite{Paris}, let $\mathcal{T}_p(\Sigma_g, \mathcal{P})$ be the kernel of the  action of the $n$-th punctured mapping class group
of  $\Sigma_g$ on $H_1(\Sigma, \F_p)$. The group  $\mathcal{T}_p(\Sigma_g, \mathcal{P})$ 
is  residually $p$-finite~\cite{Paris}.
Since $P_n(\Sigma_g)$ can be easily realised as subgroups of $\mathcal{T}_p(\Sigma_g, \mathcal{P})$,
one derives another proof of the residually $p$-finiteness of $P_n(\Sigma_g)$.

Finally, let us remark that \reth{main}  could be proved using the group presentation of $P_n(\Sigma_g)$
and the algebraic section proposed in~\cite{GG}, but related computations would become much more involved.

\section{Braid groups on $2$ strands}

Let us   recall a well known on the braid group $B_n$.

\begin{prop} \label{braidlcs}
Let $B_n$ be the Artin braid group on $n\ge 3$ strands.\\
Then $\Gamma_1(B_n)/\Gamma_2(B_n) \cong \Z$ and $\Gamma_2(B_n)=\Gamma_3(B_n)$.
\end{prop}

A similar result holds for Artin-Tits group of finite type~\cite{BGG}.

Now, let $\Sigma_g$ be a closed oriented surface of genus $g>0$ and let $B_n(\Sigma_g)$ be the braid group on $n$ strands of $\Sigma_g$.
The main result has been to determine all  lower central quotients of surface braid groups on at least $3$ strands.
In particular,  it was proven that $\Gamma_2(B_n(\Sigma_g))/\Gamma_3(B_n(\Sigma_g)) \simeq \Z_{n-1+g}$
and  that $\Gamma_3(B_n(\Sigma_g))=\Gamma_4(B_n(\Sigma_g))$ for $n \ge 3$.

Since for $n=1$ we have that $B_1(\Sigma)=\pi_1(\Sigma)$, which is residually free, in order to complete the study of lower central series of braid groups of closed surfaces we need to consider the case of $2$ strands.
In the case of the torus $\T^2$, in~\cite{BGG} was proved that $B_2(\T^2)$ is residually nilpotent using essentially the fact that this group is a central extension of $\Z_2 \ast \Z_2 \ast \Z_2$.

 Actually the residually nilpotence of $B_2(\Sigma_g)$, for any $g>0$, is a consequence of
 \reth{main} and of the following result of  Gruenberg.
 As above we set $\mathcal{FP}$ the class of groups
having the group-theoretic property $\mathcal{P}$.

 \begin{lem}{\bf (\cite{Gr})}\label{gru}
Let  $\mathcal{FP}$ be one of following classes:
\begin{enumerate}[(i)]
\item the class of solvable groups;
\item the class of finite groups;
\item the class of p-finite groups for a given prime number $p$.
\end{enumerate}
Let $P\in \mathcal{FP}$ and suppose that $H$ is residually $\mathcal{P}$.
Then for each extension $1 \to H \to G \to P \to 1$, the group $G$
is residually $\mathcal{P}$.
 \end{lem}

 \begin{cor} \label{b2}
 The group $B_2(\Sigma_g)$ is residually $2$-finite. In particular $B_2(\Sigma_g)$ is residually nilpotent.
 \end{cor}

 \begin{proof}
 Let $S_n$ the symmetric group on $n$ elements.
 We recall that $P_n(\Sigma_g)$ can be realized as the kernel of the canonical projection
 $\pi: B_n(\Sigma_g) \to S_n$.
 Since $S_2=\Z_2$ and $P_2(\Sigma_g)$ is residually torsion free, the hypothesis of Lemma~\ref{gru}
 are fulfilled and the claim follows.
 \end{proof}

We recall that pure braid groups of surfaces with non empty boundary are
 residually torsion free nilpotent: then, from Lemma~\ref{gru} it follows
 that also braid groups on $2$ strands of surfaces with non empty boundary are
 residually 2-finite.

 Therefore,  for any oriented and connected  surface $\Sigma$ of positive genus, $B_2(\Sigma)$
 is residually nilpotent.
 In~\cite{BGG} was also proved that $B_2(\T^2)$ is not residually torsion free nilpotent neither bi-orderable.
 We don't know if $B_2(\Sigma)$ is bi-orderable when $\Sigma$
 is an oriented and connected  surface of positive genus different from $\T^2$.


\begin{thebibliography}{MKS}




{\small
\bibitem[BB]{bb1}
 P.~Bellingeri and  V.~Bardakov, Combinatorial properties of virtual braids,
 math.GT/0609563.


\bibitem[B]{B} P.~Bellingeri, On presentations of surface braid
groups, \emph{J.~Algebra}, \textbf{274} (2004), 543--563.


\bibitem[BGG]{BGG} P.~Bellingeri, J.~Guaschi and S.~Gervais,
Lower central series for surface braids, math.GT/0512155.

\bibitem[CCP]{CCP} D.~Cohen, F.~Cohen and S.~Prassidis, Centralizers of Lie algebras associated to the descending
lower central series of certain poly-free groups, math.GT/0603470.

\bibitem[FR1]{FR} M.~Falk and R.~Randell,
The lower central series of a fiber type arrangement
\emph{Invent. Math.},
\textbf{82} (1985), 77--88.


\bibitem[FR2]{FR2} M.~Falk and R.~Randell,
Pure braid groups and products of free groups.
\emph{Contemp.\ Math.}, \textbf{78}(1988), 217--228.


\bibitem[F]{F}
R.~H.~Fox
Free differential calculus I: Derivation in the Free Group Ring,
\emph{Ann. of Math.}, \textbf{57}(1953), 547-560.


\bibitem[GG]{GG} D.~Gon\c{c}alves  and J.~Guaschi, The roots of the full
twist for surface braid groups, \emph{Math.\ Proc.\ Camb.\ Phil.\
Soc.} \textbf{137} (2004), 307--320.

\bibitem[Go]{Go} J.~Gonz\'alez-Meneses, Ordering pure braid groups on closed surfaces,  \emph{Pacific
J.~Math.}
\textbf{203} (2002),  n\textsuperscript{o}~2, 369--378.


 \bibitem[GP]{GP}
 J.~Gonz\'alez-Meneses et  L.~Paris,
\newblock {\em Vassiliev Invariants for braids on
 surfaces},\emph{Trans. A.M.S.}  \textbf{356} n\textsuperscript{o}~1 (2004), 219--243.



\bibitem[Gr]{Gr} K.~W.~Gruenberg, Residual properties of infinite
soluble groups, \emph{Proc.\ London Math.\ Soc.} \textbf{7} (1957),
29--62.


\bibitem[MKS]{KMS}
W. Magnus, A. Karrass, D. Solitar, Combinatorial group theory,
Interscience Publishers, New York, 1996.



\bibitem[MR]{MR} R.~Mura and A.~Rhemtulla,
 Orderable groups, Lecture Notes in Pure and Applied Mathematics \textbf{27},  Marcel Dekker, New York, 1977.

\bibitem[P]{pap} S.~Papadima,
The universal finite-type invariant for braids, with integer
coefficients, {\em Topology Appl.}, \textbf{118} (2002), 169-185.

\bibitem[Pa]{Paris}
L.~Paris,
Residually $p$ properties of mapping class groups and surface groups,
math.GT/0703703.

\bibitem[Pas]{pas}  I.~Bir S.~Passi, Group rings and Their Augmentation Ideals,
Lecture Notes In Mathematics \textbf{715}, Springer Veralg, Berlin Heidelberg New York, 1979.


}
\end{thebibliography}
\end{document}